\documentstyle{article}

\newtheorem{thm}{Theorem}[section]

\newtheorem{cor}[thm]{Corollary}

\newtheorem{lem}[thm]{Lemma}

\newtheorem{prop}[thm]{Proposition}

\newtheorem{rem}[thm]{Remark}

\newtheorem{ack}{Acknowledgment}

\newcommand{\be}{\begin{equation}}
\newcommand{\ee}{\end{equation}}
\newcommand{\bea}{\begin{eqnarray}}
\newcommand{\eea}{\end{eqnarray}}
\newcommand{\ba}{\begin{eqnarray*}}
\newcommand{\ea}{\end{eqnarray*}}
\newcommand{\bl}{\begin{lem}}
\newcommand{\el}{\end{lem}}
\newcommand{\bt}{\begin{thm}}
\newcommand{\et}{\end{thm}}
\newcommand{\bc}{\begin{cor}}
\newcommand{\ec}{\end{cor}}
\newcommand{\bp}{\begin{prop}}
\newcommand{\ep}{\end{prop}}
\newcommand{\p}{{\em Proof}. }

\newcommand{\Z}{{\bf Z}}
\newcommand{\R}{{\bf R}}
\newcommand{\C}{{\bf C}}

\newcommand{\F}{{\cal F}}
\newcommand{\ar}{\rightarrow}

\newcommand{\codim}{\mbox{\rm codim}\;}
\newcommand{\im}{\mbox{\rm im}\;}

\newcommand{\id}{\mbox{\rm id}}
\newcommand{\cl}{\mbox{\rm cl}}

\newcommand{\cinf}{C^\infty}
\newcommand{\D}{{\Delta}}
\newcommand{\wedTM}{{\bigwedge TM^\ast}}
\newcommand{\wedTF}{{\bigwedge T\F^\ast}}
\newcommand{\wedTFperp}{\bigwedge T\F^{\perp\ast}}

\newcommand{\XU}{{{\cal X}(U)}}
\newcommand{\FU}{{\F_U}}
\newcommand{\XFU}{{{\cal X}\left(\FU\right)}}
\newcommand{\XUF}{{{\cal X}\left(U,\FU\right)}}

\title{Long time behavior of leafwise heat flow for Riemannian foliations}

\author{Jes\'us A. \'ALVAREZ L\'OPEZ\thanks{Partially
supported by Xunta de Galicia, grant XUGA20701B95}\\
\footnotesize Departamento de Xeometr\'{\i}a e Topolox\'{\i}a,\\
\footnotesize Facultade de Matem\'aticas,\\
\footnotesize Universidade de Santiago de Compostela,\\
\footnotesize 15706 Santiago de Compostela, Spain
\and
Yuri A. KORDYUKOV\\
\footnotesize Department of Mathematics,\\
\footnotesize Ufa State Aviation Technical University,\\
\footnotesize 12~K.~Marx str., 450025 Ufa, Russia}

\date{}

\begin{document}

\bibliographystyle{plain}

\maketitle

\begin{abstract}
For any Riemannian foliation $\F$ on a closed manifold $M$ with an arbitrary
bundle-like metric, leafwise heat flow of differential forms
is proved to
preserve smoothness on $M$ at infinite time. This result
and its proof have consequences about the space of bundle-like metrics on
$M$, about the
dimension of the space of leafwise harmonic forms, and mainly about the
second term
of the differentiable spectral sequence of $\F$.
\end{abstract}

\section{Introduction and main results} \label{sec:intro}

For a smooth foliation $\F$ on a closed Riemannian manifold $M$,
{\em leafwise heat flow} means the evolution of an initial temperature
distribution on $M$ when the leaves are thermally isolated from each other.
With more
generality, we can consider leafwise heat flow induced by leafwise elliptic
differential operators with symmetric leading symbol. In this
paper, these operators are induced by leafwise elliptic differential
complexes in the usual way.

Let $E$ be a $\Z$-graded Riemannian vector
bundle over $M$,
and $d$ a first order leafwise elliptic differential
complex on $\cinf(E)$---the space of smooth sections of $E$. Let
$\delta$ denote the formal adjoint of $d$ on
$M$---it need not be formal adjoint of $d$ on the leaves. Then $D=d+\delta$
and $\D=D^2=d\delta+\delta d$ are symmetric differential operators on
$\cinf(E)$, and
thus essentially selfadjoint in the $L^2$-completion $L^2(E)$ of $\cinf(E)$
(Theorem~2.2 in \cite{Chernoff}). Then the spectral theorem defines the
leafwise heat
operator $e^{-t\D}$ on $L^2(E)$ for each $t\geq 0$.

The heat operator has a nice behavior on smooth sections at finite time:
According to
the work of J.~Roe in \cite{Roe87}, $e^{-t\D}$ defines a continuous
operator on $\cinf(E)$ which depends continuously on $t\in[0,\infty)$. But
the main
objective of this paper is to study the behavior of $e^{-t\D}$ on
$\cinf(E)$ at the
limit when $t\ar\infty$.

Recall that, by the spectral theorem,
$e^{-t\D}$ on $L^2(E)$ strongly converges to the orthogonal projection
$\Pi$ onto the kernel of the unbounded operator defined by $\D$ in $L^2(E)$.
By setting $e^{-\infty\D}=\Pi$, we get a continuous map
$[0,\infty]\times L^2(E)\ar L^2(E)$ given by $(t,\phi)\mapsto
e^{-t\D}\phi$, where $[0,\infty]$ is the one point compactification of
$[0,\infty)$---see Section~\ref{sec:general}.
Observe that $\Pi$ need not preserve $\cinf(E)$ because $D$ and $\Delta$
may not be
elliptic on $M$; they are only leafwise elliptic. The following result, which
is proved in Section~\ref{sec:general}, gives general conditions for $\Pi$ to
preserve $\cinf(E)$.

\bt\label{t:general}
With the above notation, suppose that, on $\cinf(E)$, there is a transversely
elliptic first order differential operator $A$, and there are zero order
differential
operators
$G$, $H$, $K$ and $L$ such that
\begin{equation}\label{e:key}
Ad\pm dA=Gd+dH\;,\quad A\delta\pm\delta A=K\delta+\delta L\;.
\end{equation}
Then $\Pi$ defines a continuous operator on $\cinf(E)$, we have
the leafwise Hodge decomposition
$$
\cinf(E)=\ker\D\oplus\overline{\im\D}
=(\ker d\cap\ker\delta)\oplus\overline{\im d}\oplus\overline{\im\delta}\;,$$
and $(t,\phi)\mapsto e^{-t\D}\phi$ defines a continuous map
$[0,\infty]\times\cinf(E)\ar\cinf(E)$.
\et

In this paper, Theorem~\ref{t:general} is only applied to {\em Riemannian
foliations}. Such
foliations are characterized by having isometric holonomy for some metric
on local
transversals. This property is equivalent to the existence of a {\em
bundle-like metric} on $M$---a metric so that the foliation is locally
defined by Riemannian submersions. Let us mention that Riemannian
foliations were introduced by
B.~Reinhart
\cite{Reinhart59}, and certain description of their structure was given by
P.~Molino
\cite{Molino82}, \cite{Molino88}.

We mainly apply Theorem~\ref{t:general} to
the leafwise differential complex constructed as follows. Let $\Omega(M)$
(or simply
$\Omega$) be the de~Rham algebra of $M$. Consider the bigrading of $\Omega$
given by
$$
\Omega^{u,v}=
\cinf\left(\bigwedge^uT\F^{\perp\ast}\otimes\bigwedge^vT\F^\ast\right)\;,\quad
u,v\in\Z\;.
$$
The de~Rham derivative and coderivative decompose as sum of
bihomogeneous components,
$$d=d_{0,1}+d_{1,0} +d_{2,-1}\;,\quad
\delta=\delta_{0,-1} +\delta_{-1,0} +\delta_{-2,1} \;,$$
where the double subindixes denote the cor\-re\-spond\-ing bidegrees. See
\cite{Alv1},
\cite{Alv2} and \cite{Alv7} for the properties of these components; in
particular
\begin{itemize}
\item each $\delta_{i,j}$ is the formal adjoint of $d_{-i,-j}$ on $M$,
\item $d_{2,-1}$ and $\delta_{-2,1}$ are of order zero,
\item $D_0=d_{0,1}+\delta_{0,-1}$ and $\Delta_0 =D_0^2$ are leafwise
elliptic and
symmetric, and
\item $D_\perp=d_{1,0}+\delta_{-1,0}$ is transversely elliptic and
symmetric.
\end{itemize}
Let
$L^2\Omega$ be the $L^2$-completion of $\Omega$, and $\Pi$ the orthogonal
projection
onto the kernel of the unbounded operator
defined by $\D_0$ in $L^2\Omega$. We prove in Section~\ref{sec:de Rham}
(Proposition~\ref{p:de Rham}) that, if $\F$ is a
Riemannian
foliation and the metric bundle-like, $(\Omega,d_{0,1})$
satisfies the conditions of Theorem~\ref{t:general} with
$D_\perp$ as $A$, obtaining the
following result which solves affirmatively a conjecture of the first author
and
P.~Tondeur \cite{AlvTond1}.

\bt\label{t:de Rham}
Let $\F$ be a Riemannian foliation on a
closed manifold $M$ with a bundle-like metric. Then
$\Pi$ defines a continuous operator on $\Omega$, we
have the leafwise Hodge decomposition
\be\label{e:leafwise Hodge}
\Omega=\ker\D_0\oplus\overline{\im\D_0}
=(\ker d_{0,1}\cap\ker\delta_{0,-1})\oplus\overline{\im d_{0,1}}
\oplus\overline{\im\delta_{0,-1}}\;,
\ee
and $(t,\alpha)\mapsto e^{-t\D_0}\alpha$ defines a continuous map
$[0,\infty]\times\Omega\ar\Omega$.
\et

Theorem~\ref{t:de Rham} is rather surprising. Indeed
we can modify a given metric by introducing a huge
bump of curvature on the plaques of a small foliation chart. Thus geodesics
in close
leaves have a sharply different behavior when going through that chart, and
so is the
behavior of leafwise heat flow on close leaves. Then one could wrongly
think this
phenomenon would break continuity of leafwise heat flow at infinite time.

The above bigrading is useful
to understand the differentiable spectral sequence $(E_i,d_i)$ of $\F$ (see
e.g. \cite{Alv1}, \cite{Alv2}): There are canonical identities
\be\label{e:E1,E2} (E_0,d_0)\equiv(\Omega,d_{0,1})\;,\quad
(E_1,d_1)\equiv(H(\Omega,d_{0,1}),d_{1,0\star})\;. \ee
Moreover the $\cinf$ topology induces a topology on each $E_i$ so that $d_i$
is
continuous. Such topology on $E_1$ need not be Hausdorff
\cite{Haefliger80}, obtaining
the bigraded subcomplex $\bar{0}_1\subset E_1$ defined as the closure of
the trivial
subspace. The quotient bigraded complex $E_1/\bar{0}_1$ will be denoted by
${\cal E}_1$.
We get
$${\cal E}_1^{u,v}\cong\Omega^{u,v}\cap\ker\D_0\;,\quad u,v\in\Z\;,$$
as direct consequence of Theorem~\ref{t:de Rham},
yielding the following dualities where $p=\dim\F$ and $q=\codim\F$: If $M$ is
oriented, then  ${\cal E}_1^{u,v}\cong{\cal E}_1^{q-u,p-v}$; if $\F$ is
oriented, then
${\cal E}_1^{u,v}\cong{\cal E}_1^{u,p-v}$;
and if $\F$ is transversely oriented, then
${\cal E}_1^{u,v}\cong{\cal E}_1^{q-u,v}$.
These isomorphisms are respectively induced by the Hodge star operators on
$\wedTM$,
$\wedTF$ and  $\wedTFperp$.

This fits into a more general setting. First observe that
$(\Omega^{0,\cdot},d_{0,1})$ can be canonically identified with the leafwise
de~Rham complex $(\Omega(\F),d_\F)$, where $\Omega (\F)=\cinf(\wedTF)$
and $d_\F$ is defined by the de~Rham derivative on the leaves. Thus
$E_1^{0,\cdot}$
can be identified with the leafwise cohomology and ${\cal E}_1^{0,\cdot}$
with the leafwise reduced cohomology. Furthermore the whole
$(\Omega,d_{0,1})$ can be considered, up to sign, as the leafwise de Rham
complex of
$\F$ with coefficients in the vector bundle $\bigwedge (TM/T\F)^\ast$ with
the flat
$\F$-partial Bott connection; i.e. the partial connection induced by the
partial Bott
connection on the normal bundle \cite{Bott71}, \cite{Bott72}. The condition
that $\F$
be Riemannian is equivalent to the existence of a metric on the normal
bundle such that
the partial Bott connection is Riemannian. Thus, with more generality, we
can consider
the leafwise de~Rham complex $(\Omega(\F,V),d_\F)$ with coefficients in any
Riemannian
vector bundle $V$ with a flat Riemannian $\F$-partial connection. The
leafwise reduced
cohomology with coefficients in $V$ will be denoted by ${\cal H}(\F,V)$.
The operator
$\delta_\F$ on $\Omega(\F,V)$, defined by the de~Rham coderivative on the
leaves, is
adjoint of $d_\F$. Let $D_\F=d_\F+\delta_\F$
and $\D_\F=D_\F^2$. Let also $L^2\Omega(\F,V)$ be the $L^2$-completion of
$\Omega(\F,V)$, and $\Pi$ the orthogonal projection of $L^2\Omega(\F,V)$
onto the kernel of $\D_\F$ in $L^2\Omega(\F,V)$.  The
following result easily follows from the case with coefficients in $\R$
(Section~\ref{sec:coefficients}).

\bc\label{c:de Rham}
Let $\F$ be a Riemannian foliation on a closed manifold $M$, and let $V$ be any
Riemannian vector bundle with a flat Riemannian $\F$-partial connection. Fix
any Riemannian metric on the leaves, smooth on $M$. Then $\Pi$
defines a continuous operator on $\Omega(\F,V)$, we have the leafwise Hodge
decomposition
$$\Omega(\F,V)=\ker\D_\F\oplus\overline{\im\D_\F} =(\ker
d_\F\cap\ker\delta_\F)\oplus\overline{\im d_\F}\oplus
\overline{\im\delta_\F}\;,$$
and $(t,\alpha)\mapsto e^{-t\D_\F}\alpha$ defines a continuous map
$[0,\infty]\times\Omega(\F,V)\ar\Omega(\F,V)$.
Thus ${\cal H}(\F,V)$ can be canonically identified with $\ker\D_\F$, and,
if $\F$ is
oriented, then the leafwise Hodge star operator on $\ker\D_\F$ induces an
isomorphism
${\cal H}^v(\F,V)\cong{\cal H}^{p-v}(\F,V^\ast)$. \ec

The case of Corollary~\ref{c:de Rham} with coefficients in the normal bundle
helps to understand infinitesimal deformations of $\F$ \cite{Heitsch75}. Another
interesting application is to take coefficients in the symmetric tensor product
$S^2((TM/T\F)^\ast)$, obtaining in Section~\ref{sec:metrics} the
following consequence which solves a problem proposed by
E.~Mac\'{\i}as\footnote{Problem~4.12
of ``Open Problems'' in: {\em Analysis
and Geometry in Foliated Manifolds}. X.~Masa, E.~Mac\'{\i}as Virg\'os and
J.~A.~\'Alvarez L\'opez (Editors). Proceedings of the VII International
Colloquium on Differential Geometry, Santiago de Compostela, 26--30 July,
1994.
World Scientific, Singapore, 1995.}.

\bc\label{c:bundle-like metrics}
Let $\F$ be a Riemannian foliation on a closed manifold $M$. Then, with
respect to
the $\cinf$ topology, the space of bundle-like metrics on $M$ is a
deformation retract
of the space of all metrics on $M$.
\ec

Corollary~\ref{c:de Rham} also has the following consequence
which partially generalizes results from \cite{AlvHurder}.

\bc\label{c:harmonic forms}
Let $\F$ be a Riemannian foliation on a closed manifold $M$, and let $V$ any
Riemannian vector bundle with a flat Riemannian $\F$-partial connection. For any
metric on the leaves, smooth on $M$, if there is a nontrivial integrable
harmonic
$r$-form on some leaf with coefficients in $V$, then
$\Omega^r(\F,V)\cap\ker\D_\F$ is of infinite dimension.  \ec

As suggested in \cite{AlvHurder}, Corollary~\ref{c:harmonic forms} may be
used to find
examples of Riemannian foliations on closed Riemannian manifolds with dense
leaves and
infinite dimensional space of smooth leafwise harmonic forms.

Theorem~\ref{t:de Rham} and its proof are also useful to get a better
understanding of
the term $E_2$ in the spectral sequence of any Riemannian foliation $\F$ on
a closed
manifold $M$. Let ${\cal H}_1=\ker\D_0$ and $\widetilde{\cal
H}_1=\overline{\im\D_0}$ in
$\Omega$, and let $L^2{\cal H}_1$ and  $L^2\widetilde{\cal H}_1$ be the
corresponding
closures in $L^2\Omega$. Consider the bigrading of ${\cal H}_1$ induced by
the one of
$\Omega$. If $\F$ is Riemannian and the metric bundle-like, by
(\ref{e:E1,E2}) and
Theorem~\ref{t:de Rham}, the differential map $d_1$ on ${\cal E}_1$ corresponds
to the map $\Pi d_{1,0}$ on ${\cal H}_1$, which will be also denoted by
$d_1$. Hence
$H^u({\cal E}_1^{\cdot,v},d_1)\cong H^u({\cal H}_1^{\cdot,v},d_1)$.
Consider also the
following operators on ${\cal H}_1$: $\delta_1=\Pi\delta_{-1,0}$,
$D_1=d_1+\delta_1$ and
$\D_1=D_1^2$. Such a $\delta_1$ is adjoint of $d_1$, and thus $D_1$ and
$\D_1$ are
symmetric in $L^2{\cal H}_1$.

We also define a differential map $\tilde d_1$ on  $\widetilde{\cal
H}_1$ as follows. First we slightly change the bigrading that
$\Omega^{\cdot,\cdot}$
induces on $\widetilde{\cal H}_1$: Set
$$\widetilde{\cal
H}_1^{u,v}=\overline{d_{0,1}(\Omega^{u,v-1})}\oplus
\overline{\delta_{0,-1}(\Omega^{u+1,v})}\;,\quad u,v\in\Z\;.$$
Let $\widetilde\Pi_{\cdot,v}$ be the
projection of $\Omega$ onto $\widetilde{\cal H}_1^{\cdot,v}$ according to
(\ref{e:leafwise Hodge}), and set $\tilde d_1= \widetilde\Pi_{\cdot,v}d$ on
$\widetilde{\cal H}_1^{\cdot,v}$. We shall see that  $\tilde d_1^2=0$
(Section~\ref{sec:E2}); indeed $H^u(\bar{0}_1^{\cdot,v},d_1)\cong
H^u(\widetilde{\cal
H}_1^{\cdot,v},\tilde d_1)$---see Section~\ref{sec:E2} for a better
understanding of
this modified bigrading. Then set
$\tilde\delta_1=\widetilde\Pi_{\cdot,v}\delta$ on
$\widetilde{\cal H}_1^{\cdot,v}$ for each $v$, and let $\widetilde D_1=\tilde
d_1+\tilde\delta_1$ and
$\widetilde\D_1=\widetilde D_1^2$. Such $\tilde\delta_1$ is adjoint of
$\tilde d_1$, and thus $\widetilde D_1$ and $\widetilde\D_1$ are symmetric in
$L^2\widetilde{\cal H}_1$. By using Theorem~\ref{t:de Rham} and the role
played by
$D_\perp$ in its proof, we prove in Section~\ref{sec:E2} the following
result which
generalizes the basic Hodge decompositions of \cite{KacimiHector} and
\cite{KamberTond87}.

\bt\label{t:E2}
Let $\F$ be a Riemannian foliation on a closed manifold $M$. For any bundle-like
metric on $M$, the operators $\D_1$ and $\widetilde\D_1$ are essentially
self-adjoint on $L^2{\cal H}_1$ and $L^2\widetilde{\cal H}_1$, respectively.
Moreover $L^2{\cal H}_1$ and  $L^2\widetilde{\cal H}_1$ have complete
orthonormal
systems, $\{\phi_i : i=1,2,\ldots\}\subset{\cal H}_1$ and
$\{\tilde\phi_i : i=1,2,\ldots\}\subset\widetilde{\cal H}_1$,
consisting of eigenvectors of $\D_1$ and $\widetilde\D_1$, respectively, so that
the corresponding eigenvalues satisfy
$0\leq\lambda_1\leq\lambda_2\leq\cdots$,
$0<\tilde\lambda_1\leq\tilde\lambda_2\leq\cdots$, with $\lambda_i\uparrow\infty$
if $\dim{\cal H}_1=\infty$, and $\tilde\lambda_i\uparrow\infty$
if $\dim\widetilde{\cal H}_1=\infty$; thus all of these eigenvalues have finite
multiplicity. In particular we have $${\cal H}_1=\ker\D_1\oplus\im\D_1\;,\quad
\widetilde{\cal H}_1=\im\widetilde\D_1\;.$$
\et

A direct consequence of Theorem~\ref{t:E2} is another proof of the
following known result, which has important implications about tautness of
Riemannian foliations \cite{Masa92}, \cite{Alv7}.

\bc[X.~Masa \cite{Masa92}; see also \cite{Alv1}, \cite{Alv2}]\label{c:E2}
Under the same conditions we have $H(\bar 0_1)=0$ and $E_2\cong H({\cal E}_1)$,
which is of finite dimension and, if $M$ is orientable, satisfies the duality
$E_2^{u,v}\cong E_2^{q-u,p-v}$, $u,v\in\Z$, where $p=\dim\F$ and $q=\codim\F$.
\ec

The duality in Corollary~\ref{c:E2} is induced by the Hodge star operator
of $M$ acting on
$\ker\D_1$ for any bundle-like metric.

Observe that Corollary~\ref{c:E2} is not satisfied by arbitrary
foliations on closed manifolds \cite{Schwarz}. So, if there is a version
of Theorem~\ref{t:de Rham} for more general foliations with the same kind
of arguments,
then $D_\perp$ should be replaced in its proof by other transversely elliptic
operator and perhaps  more general conditions should be used in
Theorem~\ref{t:general}
(see Remarks~\ref{rem:more general} and~\ref{rem:K=L}).
There are related results
for non-Riemannian foliations with very
different proofs \cite{Ledrappier}, \cite{Kanai}.

In possible generalizations of Theorem~\ref{t:de Rham}, a key role may be
played by the fact that our leafwise elliptic operators are symmetric on
$M$ instead of being symmetric on the leaves. For Riemannian foliations and
bundle-like metrics both points of view are the same. In general, the
Laplacian on the leaves acting on functions induces the physical leafwise
heat flow, while $\D_0$ induces a modification of the physical leafwise
heat flow by ``exterior influence''. For non-Riemannian foliations, the
physical leafwise heat flow may ``break'' continuous functions at infinite
time, as can be easily seen for foliations on the two dimensional torus
with several Reeb components. We hope the modified leafwise heat flow
induced by $\D_0$ has a better behavior at infinite time for more general
foliations.

We hope Theorems~\ref{t:de Rham} and~\ref{t:E2} will be useful to study
relations
between spectral sequences of Riemannian
foliations and adiabatic limits; precisely, to prove the results in Section~5 of
\cite{Forman} without the strong hypothesis that the
positive spectrum of $\D_0$ is bounded away from zero.

We fix the following notation to be used in next sections.  The $k$th Sobolev
completion of $\cinf(E)$, $\Omega$ and $\Omega(\F,V)$ will be respectively
denoted by
$W^k(E)$, $W^k\Omega$ and $W^k\Omega(\F,V)$. Fix a norm $\|\cdot\|_k$ in
any of these
$k$th Sobolev spaces, and let $\|\cdot\|_k$ also denote the corresponding
norm of
bounded operators on that space. Finally, closure in $k$th Sobolev spaces
will be denoted by $\cl_k$.

\begin{ack}
The first author would like to thank S.~Hurder and J.~Roe for helpful
conversations
about leafwise heat flow.
The second author gratefully acknowledges the hospitality and
support of the
University of Santiago de Compostela.
\end{ack}

\section{The general result}\label{sec:general}
This section is devoted to the proof of Theorem~\ref{t:general}. To do it,
a result of
J.~Roe in \cite{Roe87} is needed as a first step. We firstly state it in our setting.
Let $\cal A$ be the Fr\'echet algebra of those functions $f$ on $\R$ that
extend to
entire functions on
$\C$ such that for each compact subset $K\subset\R$ the functions
$\{x\mapsto f(x+iy)\
:\ y\in K\}$ form a bounded subset of the Schwartz space ${\cal S}(\R)$.
Such an $\cal
A$ is a module over the polynomial ring $\C[z]$, and contains all functions with
compactly supported Fourier transform and the Gaussians $x\mapsto
e^{-tx^2}$.  With the notation of Theorem~\ref{t:general}, without assuming
(\ref{e:key}), the same arguments as in Propositions~1.4 and~4.1 in
\cite{Roe87} give
the following (see also \cite{Kordyukov95} for a discussion of the action
of functions
of tangentially elliptic operators in Sobolev spaces on the ambient manifold).

\bp[{J.~Roe \cite{Roe87}}]\label{p:finite time}
The functional calculus map $f\mapsto f(D)$ given by the spectral theorem
restricts
to a continuous homomorphism of algebras and $\C[z]$-modules from $\cal A$
to the space
of bounded endomorphisms of each
$W^k(E)$, and thus to the space of continuous endomorphisms of $\cinf(E)$.
In particular, $e^{-t\D}$ defines a bounded operator on each $W^k(E)$ and a
continuous operator on $\cinf(E)$, which depends smoothly on
$t\in[0,\infty)$. 
\ep

Now assume (\ref{e:key}) is  satisfied. Then Theorem~\ref{t:general}
clearly follows
from the following six properties, which are proved by induction on
$k=0,1,2,\dots$:
\begin{itemize}

\item[(i)] There exists $c_{1,k}>0$ such that the bounded operator
$e^{-t\Delta}$ on
$W^k(E)$, defined by Proposition~\ref{p:finite time}, satisfies
$$\left\|e^{-t\Delta}\right\|_k \leq c_{1,k}\quad\mbox{for all}\quad t>0\;.$$

\item[(ii)] $De^{-t\Delta}$ defines a bounded operator on $W^k(E)$ and
there exists
$c_{2,k}>0$ such that
$$\left\|De^{-t\Delta}\right\|_k \leq
\frac{c_{2,k}}{\sqrt{t}}\quad\mbox{for all}\quad
t>0\;.$$

\item[(iii)] $\Delta e^{-t\Delta}$ defines a bounded operator on $W^k(E)$
and there
exists $c_{3,k}>0$ such that
$$\left\|\Delta e^{-t\Delta}\right\|_k \leq \frac{c_{3,k}}{t}\quad\mbox{for
all}\quad
t>0\;.$$

\item[(iv)]
The operator $e^{-t\D}$ is strongly convergent in $W^k(E)$ as
$t\ar\infty$.
Moreover $(t,\phi)\mapsto e^{-t\D}\phi$ defines a
continuous map $[0,\infty]\times W^k(E)\ar W^k(E)$.

\item[(v)] We have
\ba
W^k(E)&=&\ker\left(\D\mbox{ in } W^k(E)\right)\oplus\cl_k(\im\D)\\
&=&\ker\left(D\mbox{ in } W^k(E)\right)\oplus\cl_k(\im D)\;. \ea
The corresponding projection of $W^k(E)$ onto the kernel of $\D$ in
$W^k(E)$ is obviously defined by $\Pi$.

\item[(vi)] There exists $c_{k,4}>0$ such that
$$\| d\phi\|_k + \| \delta\phi\|_k\leq c_{k,4} \| D\phi\|_k $$
for all $\phi\in\cinf(E)$. Thus
\ba
\ker\left( D\mbox{ in } W^k(E)\right)&=
&\ker\left( d\mbox{ in }W^k(E)\right)\cap\ker\left( \delta\mbox{ in }
W^k(E)\right)\;,\\
\cl_k(\im D)&=&\cl_k(\im d)\oplus\cl_k(\im\delta)\;.
\ea

\end{itemize}

For $k=0$, properties (i)--(v) follow directly from the spectral
theorem.\medskip

{\em Proof of property~$($vi$)$ for $k=0$}. Since the image of $d$ is
orthogonal to the image of $\delta$ in
$L^2(E)$, we get
$\| D\phi\|_0^2=
\| d\phi\|_0^2+\|\delta\phi\|_0^2$ for any $\phi\in\cinf(E)$. Thus
$$\| d\phi\|_0+\|\delta\phi\|_0\leq\sqrt{2}\,\| D\phi\|_0\;.\quad\Box$$

Now suppose properties~(i)--(vi) hold for a given $k=l$ and we shall prove
them for $k=l+1$.

The direct sum decompositions in properties~(v) and~(vi) for $k=l$ define
projections $P$ and $Q$ of $W^l(E)$ onto $\cl_l(\im d)$ and
$\cl_l(\im\delta)$, respectively. Thus $\id=\Pi+P+Q$.

\bl\label{l:key'}
There are bounded operators $B_1,\ldots,B_5$ on $W^l(E)$ such that
\ba
Ad\pm dA&=&B_1d+dB_2\;,\\
A\delta\pm\delta A&=&B_1\delta+\delta B_2\;,\\
{[A,\D]}&=&B_3\D + \D B_4 + D B_5 D\;.
\ea
\el

\p  We clearly have
$$
d\Pi=\Pi d=\delta\Pi=\Pi\delta=dP=P\delta=\delta Q=Qd=0\;.
$$
So
$$
d=dQ=DQ=Pd=PD\;,\quad\delta=\delta P=DP=Q\delta=QD\;.
$$
Hence (\ref{e:key}) yields the first two equalities in the statement with
$B_1=GP+KQ$
and $B_2=QH+PL$. Thus $AD\pm DA=B_1D+DB_2$, yielding the third equality by
using the
equation
$$[A,\D]=(AD\pm DA)D\mp D(AD\pm DA)\;.\quad\Box$$

\bl \label{l:Duhamel} For any bounded operator $R:W^{l+1}(E)\ar W^l(E)$, we
have the
following Duhamel type formula:
$$[ R , e^{-t\D} ]= - \int_0^t e^{-(t-s)\D} \, [ R,\D ] \, e^{-s\D} \, ds$$
as a bounded operator $W^{l+1}(E)\ar W^l(E)$.
\el

\p This follows by arguing as in the proof of the usual Duhamel's formula
(see Lemma~12.51 in \cite{CyconFroeseKirschSimon}).~$\Box$\medskip

\bl \label{l:[A,e]}
$[A,e^{-t\Delta}]$ defines a bounded operator on $W^l(E)$, and there exists
$\tilde
c_{l,1}>0$ such that
$$\|[A,e^{-t\Delta}]\|_l\leq\tilde c_{l,1}\quad\mbox{for all}\quad t>0\;.$$
\el

\p By Lemmas~\ref{l:key'} and~\ref{l:Duhamel} we have
$[A,e^{-t\D}]=I_1+I_2$, where
\ba
I_1&=&-\int_0^te^{(t-s)\D}\,B_3\D\,e^{-s\D}\,ds-
\int_0^te^{(t-s)\D}\,\D B_4\,e^{-s\D}\,ds\;,\\
I_2&=&-\int_0^te^{(t-s)\D}\,DB_5D\,e^{-s\D}\,ds\;. \ea
On the one hand,
property~(ii) for $k=l$ yields
$$\|I_2\|_l\leq c_{l,2}^2\int_0^t\frac{ds}{\sqrt{(t-s)s}}\;,$$
which is bounded independently of $t$ since this integral is easily
seen to be $\pi$.

On the other hand,
\ba
I_1&=&e^{-t\D/2}(B_3-B_4)\,e^{-t\D/2}+e^{-t\D}\,B_3+B_4\,e^{-t\D}\\
   &&\mbox{}-\int_0^{t/2}e^{-(t-s)\D}\,\D(B_3+B_4)\,e^{-s\D}\,ds\\
   &&\mbox{}-\int_{t/2}^te^{-(t-s)\D}\,(B_3-B_4)\D\,e^{-s\D}\,ds\;,
\ea
where we have decomposed the integrals in the definition of $I_1$ as sum of
integrals
on the intervals $[0,1/2]$ and $[1/2,1]$, and we have used integration by parts
with the equalities
$$\D\,e^{-s\D}=-\frac{\partial}{\partial s}e^{-s\D}\;,\quad
\D\,e^{-(t-s)\D}=\frac{\partial}{\partial s}e^{-(t-s)\D}\;.$$
Hence, by
properties~(i) and~(iii) for $k=l$ we get
\ba
\|I_1\|_l&\leq&
c_{l,1}^2\|B_3-B_4\|_l+c_{l,1}\left(\|B_3\|_l+\|B_4\|_l\right)\\
&&\mbox{}+c_{l,1}c_{l,3}\left(\|B_3+B_4\|_l\int_0^{t/2}\frac{ds}{t-s}+
\|B_3-B_4\|_l\int_{t/2}^t\frac{ds}{s}\right)\;,
\ea
which is bounded independently of $t$ since both of these integrals are
equal to
$\ln 2$.~$\Box$\medskip

\bl \label{l:ADe+-De}
$AD\,e^{-t\D}\pm D\,e^{-t\D} A$ defines a bounded operator on $W^l(E)$ and
there exists
$\tilde c_{l,2}>0$ such that
$$\left\|AD\,e^{-t\D}\pm D\,e^{-t\D}\,A\right\|_l\leq
\frac{\tilde c_{l,2}}{\sqrt{t}}\quad\mbox{for all}\quad t>0\;.$$
\el

\p By Lemma~\ref{l:key'} we have
$$AD\pm DA=B_1D+DB_2\;,$$
yielding
\ba
\lefteqn{AD\,e^{-t\D}\pm D\,e^{-t\D}\,A}\\
&=&A\,e^{-t\D/2}\,D\,e^{-t\D/2}\pm e^{-t\D/2}\,D\,e^{-t\D/2}\,A\\
&=&\left[A,e^{-t\D/2}\right]\,D\,e^{-t\D/2}\mp
e^{-t\D/2}\,D\,\left[A,e^{-t\D/2}\right]\\
&&\mbox{}+e^{-t\D/2}\,(AD\pm DA)\,e^{-t\D/2}\\
&=&\left[A,e^{-t\D/2}\right]\,D\,e^{-t\D/2}\mp
e^{-t\D/2}\,D\,\left[A,e^{-t\D/2}\right]\\
&&\mbox{}+e^{-t\D/2}\,B_1D\,e^{-t\D/2}+e^{-t\D/2}\,DB_2\,e^{-t\D/2}\;.
\ea
Thus the result follows by Lemma~\ref{l:[A,e]} and properties~(i) and~(ii) for
$k=l$.~$\Box$\medskip

Since $D$ is a leafwise elliptic operator of order one and $A$ is transversely
elliptic of order one, there exist $e_{l,1},e_{l,2}>0$ such that
\be\label{e:norm}
e_{l,1}\|\phi\|_{l+1}\leq\|\phi\|_l+\|D\phi\|_l+\|A\phi\|_l\leq
e_{l,2}\|\phi\|_{l+1}
\ee
for all $\phi\in W^{l+1}(E)$.\medskip

{\em Proof of property~$($i$)$ for $k=l+1$}.
This property follows from (\ref{e:norm}) since, for $\phi\in W^{l+1}(E)$,
we have
the following:
$$\left\| e^{-t\D}\phi\right\|_l\leq c_{l,1}\,\|\phi\|_l\leq
c_{l,1}e_{l,2}\,\|\phi\|_{l+1}\;,$$
$$\left\| D\,e^{-t\D}\phi\right\|_l=
\left\|e^{-t\D}\,D\phi\right\|_l\leq c_{l,1}\,\|D\phi\|_l\leq
c_{l,1}e_{l,2}\,\|\phi\|_{l+1}\;,$$
\ba
\left\|A\,e^{-t\D}\phi\right\|_l&\leq&\left\|\left[A,e^{-t\D}\right]\phi
\right\|_l+
\left\|e^{-t\D}\,A\phi\right\|_l\\
&\leq&\tilde c_{l,1}\,\|\phi\|_l+  c_{l,1}\,\|A\phi\|_l\\
&\leq&(\tilde c_{l,1}+c_{l,1})e_{l,2}\,\|\phi\|_{l+1}\;,
\ea
where we have used property~(i) for $k=l$, Lemma~\ref{l:[A,e]} and
(\ref{e:norm}).~$\Box$\medskip

{\em Proof of property~$($ii$)$ for $k=l+1$}.
Again this property follows from (\ref{e:norm}) since, for $\phi\in
W^{l+1}(E)$ we
have the following:
$$\left\| D\,e^{-t\D}\phi \right\|_l\leq
\frac{c_{l,2}}{\sqrt{t}}\|\phi\|_l
\leq\frac{c_{l,2}e_{l,2}}{\sqrt{t}}\,\|\phi\|_{l+1}\;,$$
$$\left\|\D \,e^{-t\D}\phi\right\|_l=\left\|D\,e^{-t\D}\,D\phi\right\|_l\leq
\frac{c_{l,2}}{\sqrt{t}}\,\|D\phi\|_l\leq
\frac{c_{l,2}e_{l,2}}{\sqrt{t}}\,\|\phi\|_{l+1}\;,$$
\ba
\left\|AD\,e^{-t\D}\phi\right\|_l&\leq&\left\|\left(AD\,e^{-t\D}\pm
D\,e^{-t\D}\,A\right)\phi\right\|_l+\left\|D\,e^{-t\D}
\,A\phi\right\|_l\\
&\leq&\frac{\tilde c_{l,2}}{\sqrt{t}}\,\|\phi\|_l+
\frac{c_{l,2}}{\sqrt{t}}\,\|A\phi\|_l\\
&\leq&\frac{(\tilde c_{l,2}+c_{l,2})e_{l,2}}{\sqrt{t}}\,\|\phi\|_{l+1}\;,
\ea
where we have used property~(ii) for $k=l$, Lemma~\ref{l:ADe+-De} and
(\ref{e:norm}).~$\Box$\medskip

{\em Proof of property~$($iii$)$ for $k=l+1$}.
This follows directly from property~(ii) for $k=l+1$ since
$\D \,e^{-t\D}=D\,e^{-t\D/2}\,D\,e^{-t\D/2}$ on $\cinf(E)$.~$\Box$\medskip

Let $\widetilde\Pi=\id-\Pi$.
Set $B=B_4\Pi-B_3\widetilde\Pi$, which is a bounded operator on $W^l(E)$.
For further
reference, we point out the following estimates which are similar to
(\ref{e:norm}):
There exist $e_{l,1}',e_{l,2}',e_{l,1}'',e_{l,2}''>0$ such
that
\be\label{e:norm'}
e_{l,1}'\|\phi\|_{l+1}\leq\|\phi\|_l+\| D\phi\|_l+\|(A+B)\phi\|_l\leq
e_{l,2}'\|\phi\|_{l+1}\;,
\ee
\be\label{e:norm''}
e_{l,1}''\|\phi\|_{l+1}\leq\|\phi\|_l+\| D\phi\|_l+\|(A-B_1)\phi\|_l\leq
e_{l,2}''\|\phi\|_{l+1}\;.
\ee

\bl\label{l:[A+B,e]}
The operator $\left[A+B,e^{-t\D}\right]$ strongly converges to the zero
operator on $W^l(E)$ as $t\ar\infty$. Moreover 
$(t,\phi)\mapsto\left[A+B,e^{-t\D}\right]\phi$ extends to a continuous map
$[0,\infty]\times W^l(E)\ar W^l(E)$ vanishing on $\{\infty\}\times W^l(E)$.
\el

\p Take any $a>1$ that will be fixed later.  Let us write 
$\left[A+B,e^{-t\D}\right]=I_1+I_2$, where
\ba
I_1 & = &\left[A+B,e^{-(t-t/a)\D}\right] e^{-t\D/a}\;,\\
I_2 & = &e^{-(t-t/a)\D} \left[A+B,e^{-t\D/a}\right]\;.
\ea

By Lemma~\ref{l:key'} and since $\D\Pi=0$, we easily get
\be\label{e:[A+B,Delta]}
[A+B,\D]=(\D(B_3+B_4)+DB_5D)\widetilde\Pi\;.
\ee
Hence, by Lemma~\ref{l:Duhamel} we have
$I_2=I_{2,1}+I_{2,2}$, where
\ba
I_{2,1}&=&-\int_0^{t/a}e^{-(t-s)\D}\,\D(B_3+B_4)\,e^{-s\D}\,\widetilde\Pi\,ds\;,\\
I_{2,2}&=&-\int_0^{t/a}e^{-(t-s)\D}\,DB_5D\,e^{-s\D}\,\widetilde\Pi\,ds\;.
\ea
Properties~(i),~(iii) and~(vi)  for $k=l$ yield the
following estimate for $\phi\in W^l(E)$:
\ba
\| I_{2,1}\phi\|_l&\leq&
c_{l,1}c_{l,3}\,\|B_3+B_4\|_l\,\left\|\widetilde\Pi\phi\right\|_l
\int_0^{t/a}\frac{ds}{t-s}\\
&\leq&c_{l,1}c_{l,3}\,\|B_3+B_4\|_l\,
\left\|\widetilde\Pi\phi\right\|_l\,\ln\frac{a}{a-1}\;.
\ea
Similarly, properties~(i),~(ii) and~(vi) for $k=l$ yield the following 
estimate for $\psi\in W^l(E)$:
\ba
\|I_{2,2}\phi\|_l&\leq&c_{l,2}^2\left\|\widetilde\Pi\phi\right\|_l
\int_0^{t/a}\frac{ds}{\sqrt{(t-s)s}}\\
&\leq&c_{l,2}^2\left\|\widetilde\Pi\phi\right\|_l
\left(\pi+\arctan\frac{2-a}{a}\right)\;, 
\ea 
Therefore $I_2$ defines a bounded operator on $W^l(E)$ whose norm can be
made arbitrarily small uniformly on $t$ by taking $a$ large enough.

To study $I_1$ we shall use the following:  
\be\label{e:[A+B,e]Pi}
\left[A+B,e^{-s\D}\right]\Pi=0\quad\mbox{for all}\quad s\geq 0\;.
\ee
This equation holds because it is obvious for $s=0$ and moreover 
\ba
\frac{d}{ds}\left(\left[A+B,e^{-s\D}\right]\Pi\right)
&=&-\left[A+B,e^{-s\D}\D\right]\Pi\\
&=&-\left[A+B,e^{-s\D}\right]\D\Pi-e^{-s\D}[A+B,\D]\Pi\;,
\ea
which vanishes since $\D\Pi=0$ and by (\ref{e:[A+B,Delta]}).
Now (\ref{e:[A+B,e]Pi}) yields
\ba
I_1 & = &\left[A+B,e^{-(t-t/a)\D}\right] \widetilde\Pi e^{-t\D/a}\\
& = & \left[A+B,e^{-(t-t/a)\D}\right] \left(e^{-t\D/a}-\Pi\right)\;.
\ea
Furthermore, by Lemma~\ref{l:[A,e]} and property~(i) for $k=l$, 
$\left[A+B,e^{-(t-t/a)\Delta}\right]$ defines a bounded operator on $W^l(E)$
whose norm is uniformly bounded on $t$. Therefore,
by property~(iv) for $k=l$, $I_1$ strongly converges to zero
in $W^l(E)$ as $t\ar\infty$ for any $a>1$, and 
$(t,\phi)\mapsto I_1\phi$ extends to a continuous map
$[0,\infty]\times W^l(E)\ar  W^l(E)$ vanishing on
$\{\infty\}\times W^l(E)$.~$\Box$\medskip

{\em Proof of property~$($iv$)$ for $k=l+1$}.
Consider the following bounded compositions:
\be\label{e:composition1}
\mbox{\begin{picture}(166,22)
\put(0,0){$W^{l+1}(E)$}
\put(43,3){\vector(1,0){30}}
\put(46,6){$e^{-t\D}$}
\put(77,0){$W^{l+1}(E)\hookrightarrow W^l(E)\;,$}
\end{picture}}
\ee
\be\label{e:composition2}
\mbox{\begin{picture}(176,17)
\put(0,0){$W^{l+1}(E)$}
\put(43,3){\vector(1,0){30}}
\put(46,6){$e^{-t\D}$}
\put(77,0){$W^{l+1}(E)$}
\put(120,3){\vector(1,0){20}}
\put(125,6){$D$}
\put(144,0){$W^l(E)\;,$}
\end{picture}}
\ee
\be\label{e:composition3}
\mbox{\begin{picture}(192,17)
\put(0,0){$W^{l+1}(E)$}
\put(43,3){\vector(1,0){30}}
\put(46,6){$e^{-t\D}$}
\put(77,0){$W^{l+1}(E)$}
\put(120,3){\vector(1,0){36}}
\put(123,6){$A+B$}
\put(160,0){$W^l(E)\;.$}
\end{picture}}
\ee
By~(\ref{e:norm'}), it is enough to prove
that the compositions~(\ref{e:composition1}),~(\ref{e:composition2})
and~(\ref{e:composition3}) are strongly convergent as $t\ar\infty$, and to
prove the continuous extension to $[0,\infty]\times  W^{l+1}(E)$ of the maps
$[0,\infty)\times W^{l+1}(E)\ar W^l(E)$ defined by these time dependent
operators. This holds for~(\ref{e:composition1}) and~(\ref{e:composition2})
by property~(iv) for $k=l$ since these operators are respectively equal to
the compositions
\begin{center}
\begin{picture}(166,43)
\put(5,26){$W^{l+1}(E)\hookrightarrow W^l(E)$}
\put(95,29){\vector(1,0){30}}
\put(98,32){$e^{-t\D}$}
\put(129,26){$W^l(E)\;,$}
\put(0,0){$W^{l+1}(E)$}
\put(43,3){\vector(1,0){20}}
\put(48,6){$D$}
\put(67,0){$W^l(E)$}
\put(100,3){\vector(1,0){30}}
\put(103,6){$e^{-t\D}$}
\put(134,0){$W^l(E)\;.$}
\end{picture}
\end{center}

On the other hand,~(\ref{e:composition3}) is the sum
of the bounded compositions 
\be\label{e:composition4}
\mbox{\begin{picture}(196,17)
\put(0,0){$W^{l+1}(E)\hookrightarrow W^l(E)$}
\put(90,3){\vector(1,0){70}}
\put(93,9){$\left[A+B,e^{-t\D}\right]$}
\put(164,0){$W^l(E)\;,$}
\end{picture}}
\ee
\be\label{e:composition5}
\mbox{\begin{picture}(182,17)
\put(0,0){$W^{l+1}(E)$}
\put(43,3){\vector(1,0){36}}
\put(46,6){$A+B$}
\put(83,0){$W^l(E)$}
\put(116,3){\vector(1,0){30}}
\put(119,6){$e^{-t\D}$}
\put(150,0){$W^l(E)\;.$}
\end{picture}}
\ee
Now Lemma~\ref{l:[A+B,e]} and property~(iv) for $k=l$ respectively
imply the strong convergence of (\ref{e:composition4}) and
(\ref{e:composition5}) as $t\ar\infty$, as well as the continuous extension to 
$[0,\infty]\times W^{l+1}(E)$ of the maps
$[0,\infty)\times W^{l+1}(E)\ar W^l(E)$ defined by both of these
time dependent operators.~$\Box$\medskip

\bc\label{c:bounded Pi}
$\Pi$ defines a bounded operator on $W^{l+1}(E)$.
\ec

\p This result is a direct consequence of properties~(i) and~(iv) for
$k=l+1$---observe that, for each $\phi\in W^{l+1}(E)$, the limit
of $e^{-t\D}\phi$ in $W^{l+1}(E)$ as $t\ar\infty$ can only be $\Pi\phi$ since it
is so in $L^2(E)$.~$\Box$\medskip

\bc\label{c:[A+B,Pi]}
$[A+B,\Pi]=0:W^{l+1}(E)\ar W^l(E)$.
\ec

\p By Lemma~\ref{l:[A+B,e]}, $\left[A+B,e^{-t\D}\right]$
strongly converges to zero
on $W^l(E)$ as $t\ar\infty$. Hence the result follows because, as
$t\ar\infty$ and for each $\phi\in W^{l+1}(E)$,
$(A+B)\,e^{-t\D}\phi$ and $e^{-t\D}(A+B)\phi$ converge in $W^l(E)$ to
$(A+B)\Pi\phi$ and $\Pi(A+B)\phi$, respectively, by
property~(iv) for $k=l+1,l$.~$\Box$\medskip

{\em Proof of property~$($v$)$ for $k=l+1$}.
By Corollary~\ref{c:bounded Pi} and property~(v) for $k=l$, $\Pi$ defines a
projection of $W^{l+1}(E)$ onto the kernel of $\D$ in $W^{l+1}(E)$.

Now, for $t>0$ and
$s\in\R$ let
$$
f_t(s)=\left\{
\begin{array}{cl}
(1-e^{-ts^2})/s^2& \mbox{if}\quad s\neq 0\\
t& \mbox{if}\quad s=0\;.
\end{array}
\right.
$$
It is easy to  check that $f_t$ is in the algebra $\cal A$ of
Proposition~\ref{p:finite time}. Thus $f_t(D)$ defines a bounded operator on
$W^{l+1}(E)$ satisfying
$$\id-e^{-t\D}=\D\,f_t(D)\;.$$
So property~(iv) for $k=l+1$ yields that,  for any $\phi\in
W^{l+1}(E)$, $\widetilde\Pi\phi$ is the limit of
$\D f_t(D)\phi$ in $W^{l+1}(E)$ as $t\ar\infty$. Hence
$$\widetilde\Pi(W^{l+1}(E))\subset\cl_{l+1}(\im\D)\;.$$
This is really an equality because the reverse inclusion can be easily proved as
follows:
\ba
\cl_{l+1}(\im\D)&\subset&W^{l+1}(E)\cap\cl_0(\im\D)\\
&=&W^{l+1}(E)\cap\widetilde\Pi\left(L^2(E)\right)\\
&=&\widetilde\Pi\left(W^{l+1}(E)\right)\;.
\ea
Here we have used property~(v) for $k=0$ and the fact
that $\widetilde\Pi$ is a projection. Therefore we have proved the
first direct sum decomposition in property~(v) for $k=l+1$.

The second direct sum decomposition follows similarly by using the
functions given by
$$
g_t(s)=\left\{
\begin{array}{cl}
(1-e^{-ts^2})/s& \mbox{if}\quad s\neq 0\\
t& \mbox{if}\quad s=0
\end{array}
\right.
$$
instead of the $f_t$.~$\Box$\medskip

{\em Proof of property~$($vi$)$ for $k=l+1$}. Take any $\phi\in\cinf(E)$.
Lemma~\ref{l:key'}, property~(vi) for $k=l$ and
(\ref{e:norm''}) yield
\ba
\|d\phi\|_{l+1}+\|\delta\phi\|_{l+1}&\leq&
\frac{1}{e_{l,1}''}\left(\|d\phi\|_l+\|\delta\phi\|_l+
\|Dd\phi\|_l+\|D\delta\phi\|_l\right.\\
& &\left.\mbox{}+\|(A-B_1)d\phi\|_l+\|(A-B_1)\delta\phi\|_l\right)\\
&\leq&\frac{1}{e_{l,1}''}\left(\|d\phi\|_l+\|\delta\phi\|_l+
\|dD\phi\|_l+\|\delta D\phi\|_l\right.\\
& &\left.\mbox{}+\|d(\mp A+B_2)\phi\|_l+\|\delta(\mp A+B_2)\phi\|_l\right)\\
&\leq&\frac{c_{l,4}}{e_{l,1}''}\left(\|D\phi\|_l+
\|\D\phi\|_l+\|D(\mp A+B_2)\phi\|_l\right)\\
&\leq&\frac{c_{l,4}}{e_{l,1}''}\left(\|D\phi\|_l+
\|\D\phi\|_l+\|(A-B_1)D\phi\|_l\right)\\
&\leq&\frac{c_{l,4}e_{l,2}''}{e_{l,1}''}
\,\|D\phi\|_{l+1}\;.\quad\Box
\ea

\begin{rem}\label{rem:symmetric A}
If $A$ is symmetric in Theorem~\ref{t:general}, then one of the two equations in
$($\ref{e:key}$)$ can be removed in its statement because both of them
are equivalent  by taking adjoints on $M$. \end{rem}

\begin{rem}\label{rem:more general}
The following more general condition can be
used instead of $($\ref{e:key}$)$ in Theorem~\ref{t:general} to get the
same long time
behavior of leafwise heat flow with a similar proof: There is some first order
transversely elliptic differential operator $A$ and zero order operators
$K_1,\ldots,K_{16}$ such that the operators
\ba F_1&=&Ad\pm dA-K_1d-dK_2-K_3\delta-\delta K_4 \;,\\
F_2&=&A\delta\pm\delta A-K_5\delta-\delta K_6-K_7d-dK_8
\ea
satisfy
\ba
\lefteqn{F_1\delta\mp\delta F_1+F_2d\mp dF_2=}\\
&&dK_9d+\delta K_{10}d+dK_{11}\delta+\delta
K_{12}\delta+d\delta K_{13}+\delta dK_{14}+K_{15}d\delta+K_{16}\delta d\;.
\ea
Nevertheless, so far we did not find any non-Riemannian
foliation with a $($non-elliptic$)$ leafwise elliptic differential complex
satisfying such a more general
condition. \end{rem}

Corollary~\ref{c:[A+B,Pi]} yields $[A,\Pi]=-[B,\Pi]$ on each $W^k(E)$,
obtaining the
following consequence that will be used later.

\bc\label{c:[A,Pi]}
$[A,\Pi]$ defines a bounded operator on each $W^k(E)$, and thus a
continuous operator
on $\cinf(E)$.
\ec

\section{Case of the leafwise de~Rham complex for Riemannian foliations}
\label{sec:de Rham}

With the notation introduced in Section~\ref{sec:intro},  the objective of
this section is to prove the following result, which
implies Theorem~\ref{t:de Rham} by using Theorem~\ref{t:general} and
Remark~\ref{rem:symmetric A}.

\bp\label{p:de Rham}
If $\F$ is a Riemannian foliation and $M$ is endowed with
a bundle-like metric, there is a zero
order differential operator $K$ on $\Omega$ such that
$$D_\perp\delta_{0,-1}+\delta_{0,-1} D_\perp
=K\delta_{0,-1}+\delta_{0,-1} K\;.$$
\ep

To prove Proposition~\ref{p:de Rham}, choose any open subset $U\subset M$ of
triviality for $\F$. Let $n=\dim M$, $p=\dim\F$ and $q=\codim\F$. Fix
tangential and
transverse orientations for $\F$ in $U$, obtaining the Hodge star operators
$\star_\F$
and $\star_\perp$ on the restrictions of $\wedTF$ and $\wedTFperp$ to $U$,
respectively. Moreover we get an induced orientation of $U$ so that
$\star_\perp(1)\wedge\star_\F(1)$ is a positive volume form. The following 
lemma can be easily proved (the statement of Lemma~4.8 in \cite{AlvTond1} is
similar).

\bl \label{l:*} Over $U$, the Hodge star operator on
$\wedTM\equiv\wedTFperp\otimes\wedTF$ is given by
$$\star=(-1)^{(q-u)v}\star_\perp\otimes\star_\F:
\bigwedge^uT\F^{\perp\ast}\otimes\bigwedge^vT\F^\ast\longrightarrow
\bigwedge^{q-u}T\F^{\perp\ast}\otimes\bigwedge^{p-v}T\F^\ast\;.$$
\el

Let $\XFU\subset\XU$ be the Lie subalgebra of vector fields which 
are tangent to the leaves of $\FU=\F|_U$, and let $\XUF\subset\XU$ be its
normalizer---the Lie algebra of infinitesimal transformations of $\FU$. 
Let also
$\Omega\left(U/\FU\right)\subset\Omega^{\cdot,0}(U)$ denote 
the basic complex of $\FU$,
and $\cinf\left(U/\FU\right)= \Omega^0\left(U/\FU\right)$. 
For $X\in\XU$ let
${\cal L}_X$ denote the corresponding Lie derivative on 
$\Omega(U)$ and $\theta_X$ its
bihomogeneous
$(0,0)$-component. By  comparing bidegrees on Cartan's 
formula we get that, if $X$ is
orthogonal to $\F$, then
$$\theta_X=i_Xd_{1,0}\quad\mbox{on}\quad\Omega^{0,\cdot}(U)\;,$$
yielding
$$\theta_{fX}=f\theta_X\quad\mbox{on}\quad\Omega^{0,\cdot}(U)\quad\mbox{for
all}\quad f\in\cinf(M)$$
and
$$d_{1,0}\beta=\sum_{i=1}^q\alpha_i\wedge\theta_{X_i}\beta\quad\mbox{for
all}\quad\beta\in\Omega^{0,\cdot}(U)\;,$$
where $X_1,\dots,X_q$ is any frame of $T\F^\perp$ on $U$ with dual coframe
$\alpha_1,\dots,\alpha_q$ in $\Omega^{1,0}(U)$. Furthermore, if $X\in\XUF$ is
orthogonal to the leaves and $Y\in\XFU$, then 
$$i_Y{\cal L}_X={\cal
L}_Xi_Y-i_{[X,Y]}=0\quad\mbox{on}\quad\Omega^{\cdot,0}(U)\;,$$ 
yielding that the
$(-1,1)$-bihomogeneous component of ${\cal L}_X$ vanishes on  
$\Omega^{\cdot,0}(U)$,
and thus on $\Omega(U)$. Therefore
$$\theta_Xd_{0,1}=d_{0,1}\theta_X\quad\mbox{on}\quad\Omega(U)$$ 
by comparing
bidegrees on the formula ${\cal L}_Xd=d{\cal L}_X$.

As was pointed out in Section~\ref{sec:intro}, the restriction 
of $\delta_{0,-1}$ to
$\Omega^{0,\cdot}\equiv\Omega(\F)$ is defined by the de~Rham 
coderivative on the
leaves. This holds whenever the transverse Riemannian volume 
element is holonomy
invariant, and in particular when the metric is bundle-like: 
On $U$,
$\omega=\star_\perp(1)$ satisfies $d\omega=0$, and thus, 
by Lemma~\ref{l:*}, for
$\beta\in\Omega^{0,v}(U)$ we have
\ba 
\delta_{0,-1}\beta&=&(-1)^{nv+n+1}\star d_{0,1}\star\beta\\
&=&(-1)^{pv+n+1}\star d_{0,1}(\omega\wedge\star_\F\beta)\\
&=&(-1)^{pv+p+1}\star(\omega\wedge d_{0,1}\star_\F\beta)\\
&=&(-1)^{pv+p+1}\star_\F d_{0,1}\star_\F\beta\;.
\ea 

\bl \label{l:[theta,delta0-1]}
If $X\in\XUF$ is orthogonal to the
leaves, then there is a zero order differential operator $R_{U,X}$ on
$\Omega^{0,\cdot}(U)$, depending $\cinf(M/\FU)$-linearly on $X$, such that 
$$[\theta_X,\delta_{0,-1}]=[R_{U,X},\delta_{0,-1}]\quad\mbox{on}\quad
\Omega^{0,\cdot}\;.$$ 
Moreover the assignment $(U,X)\mapsto R_{U,X}$ can be chosen so 
that the restriction 
of $R_{U,X}$ to any subset $U'\subset U$ of triviality for $\F$ is
equal to $R_{U',X'}$, where $X'=X|_{U'}$.
\el 

\p On $\Omega^{0,v}(U)$, since $\star_\F^2=(-1)^{(p-v)v}\id$ we get
\ba 
[\theta_X,\delta_{0,-1}]
&=&(-1)^{pv+p+1}[\theta_X,\star_\F d_{0,1}\star_\F]\\
&=&(-1)^{pv+p+1}\left([\theta_X,\star_\F]d_{0,1}\star_\F
+\star_\F d_{0,1}[\theta_X,\star_\F]\right)\\
&=&(-1)^{(p-v+1)(v-1)}[\theta_X,\star_\F]\star_\F\delta_{0,-1}\\
&&\mbox{}+(-1)^{(p-v)v}\delta_{0,-1}\star_\F[\theta_X,\star_\F]\;.
\ea 
Hence the result follows by choosing
$$R_{U,X}=(-1)^{(p-v)v}[\theta_X,\star_\F]\star_\F$$
because 
$$0=[\theta_X,\star_\F^2]=
[\theta_X,\star_\F]\star_\F+\star_\F[\theta_X,\star_\F]\;.\quad\Box$$

Obviously, $d_{0,1}$ and $\delta_{0,-1}$ are
$\cinf(U/\FU)$-linear on $\Omega(U)$. Indeed we have the following.

\bl \label{l:d01,del0-1}
We have
$$d_{0,1}\equiv(-1)^u\,\id\otimes d_{0,1}\;,\quad
\delta_{0,-1}\equiv(-1)^u\,\id\otimes\delta_{0,-1}$$
with respect to the canonical decomposition  
$$\Omega^{u,\cdot}(U)\equiv
\Omega^u\left(U/\FU\right)\otimes\Omega^{0,\cdot}(U)$$
as tensor product of $\cinf\left(U/\FU\right)$-modules.
\el 

\p The first identity is clear because $d_{0,1}$ vanishes 
on basic forms. The
second identity holds because the metric is bundle-like on $U$: 
This is equivalent to
$$\star_\perp\left(\Omega\left(U/\FU\right)\right)\subset
\Omega\left(U/\FU\right)\;,$$ 
and thus, by Lemma~\ref{l:*}, for 
$\alpha\in\Omega^u(U/\FU)$, $\beta\in\Omega^{0,v}(U)$ and $r=u+v$, we have 
\ba 
\delta_{0,-1}(\alpha\wedge\beta)&=&(-1)^{nr+n+1}\star
d_{0,1}\star(\alpha\wedge\beta)\\ &=&(-1)^{nr+n+1+(q-u)v}\star
d_{0,1}(\star_\perp\alpha\wedge\star_\F\beta)\\
&=&(-1)^{nr+n+1+(q-u)(v+1)}\star(\star_\perp\alpha\wedge 
d_{0,1}\star_\F\beta)\\
&=&(-1)^u\,\alpha\wedge\delta_{0,-1}\beta\;.\quad\Box
\ea 

The proof of Proposition~\ref{p:de Rham} can be completed as follows. Let
$X_1,\dots,X_q\in\XUF$ be a frame of $T\F^\perp$ on $U$, and let
$\alpha_1,\dots,\alpha_q\in\Omega^1\left(U/\FU\right)$ be the dual coframe. 
Take any  $\alpha\in\Omega^u\left(U/\FU\right)$ and any 
$\beta\in\Omega^{0,\cdot}(U)$.
Then Lemmas~\ref{l:[theta,delta0-1]} and~\ref{l:d01,del0-1} yield 
\ba 
(d_{1,0}\delta_{0,-1}+\delta_{0,-1}d_{1,0})(\alpha\wedge\beta)
&=&(-1)^u\,d_{1,0}\alpha\wedge\delta_{0,-1}\beta
+\alpha\wedge d_{1,0}\delta_{0,-1}\beta\\
&&\mbox{}+(-1)^{u+1}\,d_{1,0}\alpha\wedge\delta_{0,-1}\beta
-\alpha\wedge\delta_{0,-1} d_{1,0}\beta\\
&=&\sum_{i=1}^q\alpha\wedge\alpha_i\wedge[\theta_{X_i},\delta_{0,-1}]\beta\\
&=&\sum_{i=1}^q\alpha\wedge\alpha_i\wedge
[R_{U,X_i},\delta_{0,-1}]\beta\\
&=&(K_U\delta_{0,-1}+\delta_{0,-1}K_U)(\alpha\wedge\beta)\;,
\ea 
where, for $\alpha$ and $\beta$ as above,
$$
K_U(\alpha\wedge\beta)=
(-1)^u\,\sum_{i=1}^q\alpha\wedge\alpha_i\wedge R_{U,X_i}\beta\;.
$$
By the properties of the $R_{U,X_i}$, it is clear that $K_U$
is independent of the choices of the $X_i$, and moreover the restriction of 
$K_U$ to any open subset $U'\subset U$ of triviality for $\F$ is equal to
$K_{U'}$. Thus the $K_U$ defines a global operator
$K$ on $\Omega$ satisfying
$$d_{1,0}\delta_{0,-1}+\delta_{0,-1} d_{1,0}
=K\delta_{0,-1}+\delta_{0,-1} K\;,$$
which finishes the proof of 
Proposition~\ref{p:de Rham} since (see e.g. \cite{Alv1})
$$d_{1,0} d_{0,1}+d_{0,1} d_{1,0}=0\;.$$

\begin{rem}\label{rem:K=L}
Proposition~\ref{p:de Rham} is slightly stronger than~$($\ref{e:key}$)$
in Theorem~\ref{t:general}. In this case $K=L$, and thus 
$$D_\perp D_0+D_0 D_\perp= B_1D_0+D_0B_2$$
with
$$B_1=K^\ast P+KQ\;,\quad B_2=QK^\ast+PK$$
in Lemma~\ref{l:key'}. In particular $B_2=B_1^\ast$. 
Observe that, if $K$ is symmetric,
then so is $B_1$, yielding that $D_\perp-B_1$ commutes with 
$\D_0$ and thus with
$e^{-t\D_0}$. Therefore the proof of Theorem~\ref{t:de Rham} 
would follow with a much
simpler induction argument. Of course $K$ is not symmetric in general.
\end{rem}

\begin{rem}
The above proof of Proposition~\ref{p:de Rham} gives explicit
local descriptions of $K$ and $L$. But an alternative 
proof can be made by using
Molino's description of Riemannian foliations: 
It allows to reduce the proof to the
case of transversely parallelizable foliations, 
where our local arguments can be made
globally. 
\end{rem}

\begin{rem}\label{rem:[Dperp,Pi]}
In the setting of this section, Corollary~\ref{c:[A,Pi]} states that 
$[D_\perp,\Pi]$
defines a bounded operator on each $W^k\Omega$.
\end{rem}

\section{Case of leafwise differential forms with appropriate
coefficients}\label{sec:coefficients}

Corollary~\ref{c:de Rham} is proved in this section. Let thus $\F$,
$M$ and $V$ be as in the statement of that result. Since any metric 
on the leaves,
smooth on $M$, can be extended to a bundle-like metric on $M$, then 
Corollary~\ref{c:de Rham} follows directly from Theorem~\ref{t:de Rham} 
when $V$ is
any trivial Riemannian vector bundle
with the trivial $\F$-partial connection. 

In the general case we follow Molino's idea to describe Riemannian 
foliations \cite{Molino82}, \cite{Molino88}. Let $\pi:F\ar M$ be the principal
$O(k)$-bundle of orthonormal frames of $V$, where $k$ is the rank of $V$;
observe that such an $F$ is  a closed manifold. The metric $\F$-partial
connection on $V$ can be understood as an $O(k)$-invariant vector subbundle
$H\subset TF$ so that $\pi_\ast:H_f\ar T_{\pi(f)}\F$ is an isomorphism for
every frame $f\in F$. Moreover the flatness of the connection means that $H$
defines a completely integrable distribution, and let thus $\widehat\F$ be the
corresponding foliation on $F$. It is clear that $\widehat\F$ is also a
Riemannian foliation, $\pi^\ast V$ has a canonical trivialization as Riemannian
vector bundle, the pull-back of the partial connection on $V$ is the trivial
$\widehat\F$-partial connection on $\pi^\ast V$, and $\pi^\ast$ defines an
injection $\Omega(\F,V)\hookrightarrow\Omega(\widehat\F,\pi^\ast V)$. Moreover
it is easy to check that, for the lift of any given metric on the leaves of
$\F$ to the leaves of $\widehat\F$, $\D_\F$ is the restriction of
$\D_{\widehat\F}$ by the above injection. Therefore Corollary~\ref{c:de Rham}
for $\F$ and $V$ follows from the case of $\widehat\F$ and $\pi^\ast V$.

\section{The space of bundle-like metrics}\label{sec:metrics}
We prove Corollary~\ref{c:bundle-like metrics} in this section. First we
recall some technicalities from  \cite{Sanmartin95}. For a given
foliation $\F$ on a manifold $M$, let $\nu=TM/T\F$, $Q=S^2(\nu^\ast)$, and
$Q^+\subset Q$ the subbundle given by the positive definite elements in
$Q$. Hence $\cinf(Q^+)$ is the space of metrics on the normal bundle
$\nu$. Such metrics are the key point to prove 
Corollary~\ref{c:bundle-like metrics} because any metric $g$ on $M$ is
uniquely determined by fixing three objects: A metric $g_\F$ on $T\F$, a
subbundle $N\subset TM$ which is complementary of $T\F$, and a metric
$g_{\nu}$ on $\nu$. In fact, $g_{\nu}$ determines a metric $g_N$ on $N$
by the canonical isomorphism $N\cong \nu$, and $g$ is determined as the
orthogonal sum of $g_\F$ and $g_N$. Conversely, $g$ determines
$g_\F=g|_{T\F}$, $N=T\F^\perp$ and $g_{\nu}$ as the only metric that
corresponds to $g|_{N}$ by the above isomorphism. According to this
notation, the metric $g$ is bundle-like if and only if the corresponding
metric $g_{\nu}$ is parallel with respect to the $\F$-partial Bott
connection on $S^2(\nu^\ast)$; i.e. $d_\F(g_{\nu})=0$.  Thus, by
modifying only $g_{\nu}$ for each metric $g$, it is clear that
Corollary~\ref{c:bundle-like metrics} follows from the following result
by using Corollary~\ref{c:de Rham} with $V=Q$.

\bl
Suppose $\F$ is Riemannian and $M$ is closed with a fixed bundle-like metric.
Then the corresponding leafwise heat operator $e^{-t\D_\F}$ on $\Omega(\F,Q)$
preserves $\cinf(Q^+)$ for each $t\in[0,\infty]$.  \el

\p Consider the metric $\bar g$ on $\nu$ determined as above by the 
bundle-like metric on $M$. Let $\nu^1$ be the sphere bundle over $M$ given by
the normal vectors of $\bar g$-norm one. Then the result follows by checking
that, for any $g\in\cinf(Q^+)$, we have
$$\min_{v\in\nu^1}\left(e^{-t\D_\F}g\right)(v,v)\geq\min_{v\in\nu^1}g(v,v)
\quad\mbox{for all}\quad t\in[0,\infty)\;.$$
This in turn follows by checking that, for any given $T\in[0,\infty)$, 
if the map $v\in\nu^1\mapsto (e^{-T\D_\F}g)(v,v)\in\R$ reaches the minimum at
some $v_m\in\nu^1_x$ for some $x\in M$, then  
$$\left.\frac{\partial\left(e^{-t\D_\F}g\right)}{\partial t}
(v_m,v_m)\right|_{t=T}\geq 0\;.$$  
This property can be proved as follows. Extend $v_m$ to a local field $V$
of normal vectors of $\bar g$-norm one satisfying $d_\F V=0$; this is  
always possible on some open subset $U\subset M$ of triviality for $\F$ 
since $d_\F\bar g=0$. Then, if $P$ is the plaque in $U$ containing $x$, the
restriction $f_t$ of $\left(e^{-t\D_\F}g\right)(V,V)$ to $P$ satisfies the
parabolic  equation $\partial f_t/\partial t+\D_Pf_t=0$ and $f_T$  reaches the
minimum at $x$; here  $\D_P$ is the Laplacian on $P$. Hence 
$$\left.\frac{\partial f_t}{\partial t}(x)\right|_{t=T}\geq 0$$ by the maximum
principle and the proof is completed.~$\Box$\medskip

\section{Dimension of the space of leafwise harmonic forms}

Corollary~\ref{c:harmonic forms} is proved in this section. With the  notation
of that corollary, let $\phi$ be a nontrivial integrable harmonic $r$-form on
some leaf $L$ with coefficients on $V$. Such $\phi$ determines a
continuous linear functional $\tilde\phi$ on $\Omega^{p-r}(\F,V^\ast)$,
$p=\dim\F$, given by $$\tilde\phi(\psi)=\int_L\phi\wedge\psi|_{L}\;.$$
Thus $\tilde\phi$ is a singular element in $W^k\Omega^r(\F,V)$ for some
negative $k\in\Z$.  Take any sequence $\phi_i$ in $\Omega^r(\F,V)$ converging 
to $\phi$ in $W^k\Omega^r(\F,V)$. By Corollary~\ref{c:de Rham}, $\Pi\phi_i$ is
a  sequence in $\Omega^r(\F,V)\cap\ker\D_\F$ converging to the singular
$\tilde\phi=\Pi\tilde\phi$ in $W^k\Omega^r(\F,V)$, and so the $\Pi\phi_i$
generate a space of infinite dimension.

\section{Application to the second term of the spectral sequence of 
Riemannian foliations}\label{sec:E2}

The objective of this section is to prove Theorem~\ref{t:E2} and
Corollary~\ref{c:E2}. Consider thus the notation and conditions 
introduced to state
those results. 

\bl\label{l:self-adjoint}
$D_1$ and $\widetilde D_1$ are essentially self-adjoint in $L^2{\cal
H}_1$ and $L^2\widetilde{\cal H}_1$, respectively. \el

\p By Theorem~2.2 in \cite{Chernoff}, $D_\perp$ is essentially
self-adjoint in $L^2\Omega$. Then, by using e.g. Lemma~XII.1.6--(c) in
\cite{DunfordSchwartz}, so is $\Pi D_\perp\Pi$ because $\Pi$ 
is a bounded self-adjoint
operator on $L^2\Omega$. But $\Pi
D_\perp\Pi$ is equal to $D_1$ in $L^2{\cal H}_1$ 
and vanishes in its orthogonal
complement. Hence $D_1$ is essentially self-adjoint.

The proof that $\widetilde D_1$ is essentially 
self-adjoint is similar.~$\Box$\medskip

\bl\label{l:Sobolev}
We have the following properties:
\begin{itemize}

\item[$($i$)$] $D\Pi-\Pi D_\perp\Pi$ defines a bounded operator on $L^2\Omega$.

\item[$($ii$)$] For each $v\in\Z$, 
$D\widetilde\Pi_{\cdot,v}
-\widetilde\Pi_{\cdot,v}D\widetilde\Pi_{\cdot,v}$ defines a
bounded operator on $L^2\Omega$.

\end{itemize}
\el

\p Property~(i) can be proved as follows. Because $D_0\Pi=0$, we get 
$$D\Pi-\Pi D_\perp\Pi=(\id-\Pi)D_\perp\Pi+(d_{2,-1}+\delta_{-2,1})\Pi\;,$$
which is bounded on $L^2\Omega$ by Remark~\ref{rem:[Dperp,Pi]} and the formula
$$(\id-\Pi)D_\perp\Pi=[D_\perp,\Pi]\,\Pi\;.$$ 

The proof of property~(ii) is slightly more complicated. We have
$$D\widetilde\Pi_{\cdot,v}-\widetilde\Pi_{\cdot,v}D\widetilde\Pi_{\cdot,v}=
\left[D,\widetilde\Pi_{\cdot,v}\right]\widetilde\Pi_{\cdot,v}\;.$$
But, if $\phi=\phi_1+\phi_2\in\widetilde{\cal H}_1^{\cdot,v}$ with
$\phi_1\in\overline{d_{0,1}(\Omega^{\cdot,v-1})}$ and 
$\phi_2\in\overline{\delta_{0,-1}(\Omega^{\cdot,v})}$, then
$$d_{0,1}\phi_1=\delta_{0,-1}\phi_2=0\;,
\quad\delta_{0,-1}\phi_1,\ d_{0,1}\phi_2\in
\widetilde{\cal H}_1^{\cdot,v}\;,$$ 
yielding
$[D_0,\widetilde\Pi_{\cdot,v}]\phi=0$. Hence
$$\left[D,\widetilde\Pi_{\cdot,v}\right]\widetilde\Pi_{\cdot,v}=
\left[D_\perp,\widetilde\Pi_{\cdot,v}\right]\widetilde\Pi_{\cdot,v}+
\left[d_{2,-1}+\delta_{-2,1},\widetilde\Pi_{\cdot,v}\right]
\widetilde\Pi_{\cdot,v}\;.$$
But clearly
$$\left[D_\perp,\widetilde\Pi_{\cdot,v}\right]\widetilde\Pi_{\cdot,v}
=\widetilde\Pi_{\cdot,v-1}D_\perp\widetilde\Pi_{\cdot,v}
+\widetilde\Pi_{\cdot,v+1}D_\perp\widetilde\Pi_{\cdot,v}\;.$$
Therefore the result follows once we have proved that the last 
two terms define
bounded operators on $L^2\Omega$. In fact, by taking adjoints, 
it is enough to prove
that one of them defines a bounded operator for arbitrary $v$.

Let $\phi=\phi_1+\phi_2$ as above. Then obviously 
$\widetilde\Pi_{\cdot,v-1}D_\perp\phi_1=0$. 
On the other hand, $\phi_2$ is the $\cinf$
limit of $\delta_{0,-1}\psi_i$ for some sequence 
$\psi_i$ in  $\Omega^{\cdot,v}$.
So 
$$D_\perp\phi=
\lim_i\left(-\delta_{0,-1}D_\perp
+K\delta_{0,-1}+\delta_{0,-1}K\right)\psi_i$$
by
Proposition~\ref{p:de Rham}, 
yielding
$$\widetilde\Pi_{\cdot,v-1}D_\perp\phi=\widetilde\Pi_{\cdot,v-1}K\phi$$
because $\widetilde\Pi_{\cdot,v-1}\delta_{0,-1}\Omega^{\cdot,v}=0$ and
$\widetilde\Pi_{\cdot,v-1}$ is continuous on $\Omega$. Thus 
$\widetilde\Pi_{\cdot,v-1}D_\perp\widetilde\Pi_{\cdot,v}$ is bounded on
$L^2\Omega$.~$\Box$\medskip

Define the norms $\|\cdot\|_{D_1,k}$ on ${\cal H}_1$ and 
$\|\cdot\|_{\widetilde D_1,k}$ 
on  $\widetilde{\cal H}_1$ by setting
$$\|\phi\|_{D_1,k}=\left\|(\id+D_1)^k\phi\right\|_0\;,\quad
\|\psi\|_{\widetilde D_1,k}
=\left\|\left(\id+\widetilde D_1\right)^k\psi\right\|_0\;,$$
and let $W^k{\cal H}_1$ and $W^k\widetilde{\cal H}_1$ be the
corresponding completions of ${\cal H}_1$ 
and $\widetilde{\cal H}_1$. Then the
following result follows directly from Lemma~\ref{l:Sobolev}.

\bc\label{c:Sobolev}
The restrictions of the $k$th Sobolev norm 
$\|\cdot\|_k$ to ${\cal H}_1$ and 
$\widetilde{\cal H}_1$ are respectively equivalent 
to the norms $\|\cdot\|_{D_1,k}$ and
$\|\cdot\|_{\widetilde D_1,k}$. Thus $W^k{\cal H}_1$ 
and $W^k\widetilde{\cal H}_1$ are
the closures of ${\cal H}_1$ and $\widetilde{\cal H}_1$ 
in $W^k\Omega$, respectively.
\ec

The inclusions $W^{k+1}{\cal H}_1\hookrightarrow W^k{\cal H}_1$
and $W^{k+1}\widetilde{\cal H}_1\hookrightarrow 
W^k\widetilde{\cal H}_1$ are compact
operators by Corollary~\ref{c:Sobolev}. 
Then, by Proposition~2.44 in \cite{AlvTond1} and
Lemma~\ref{l:self-adjoint}, the Hilbert spaces $L^2{\cal H}_1$ and 
$L^2\widetilde{\cal H}_1$ have complete orthonormal systems, 
$\{\phi_i : i=1,2,\ldots\}\subset{\cal H}_1$ and  
$\{\tilde{\phi}_i : i=1,2,\ldots\}\subset\widetilde{\cal H}_1$,  
consisting of
eigenvectors of $\D_1$ and $\widetilde\D_1$, respectively, so that the
corresponding eigenvalues satisfy 
$0\leq\lambda_1\leq\lambda_2\leq\cdots$, 
$0\leq\tilde\lambda_1\leq\tilde\lambda_2\leq\cdots$, 
with $\lambda_i\uparrow\infty$
if $\dim{\cal H}_1=\infty$, and $\tilde\lambda_i\uparrow\infty$
if $\dim\widetilde{\cal H}_1=\infty$. Thus
it only remains to check that $\tilde\lambda_1>0$ to complete
the proof of Theorem~\ref{t:E2};
i.e. to check that $\ker\widetilde\D_1=0$. The rest of
this section will be devoted to prove this property.

For each $v\in\Z$ let
\ba
{\cal Z}_v&=&\bigoplus_{w<v}\Omega^{\cdot,w}
\oplus\ker\left(d_{0,1}:\Omega^{\cdot,v}\ar\Omega^{\cdot,v+1}\right)\;,\\
{\cal B}_v&=&\bigoplus_{w<v}\Omega^{\cdot,w}
\oplus d_{0,1}\left(\Omega^{\cdot,v-1}\right)\;.
\ea
The bigrading of $\Omega$ is used to define 
these spaces only for the sake of
simplicity, indeed they depend only on $\F$ 
\cite{Sergiescu86}, \cite{Alv2}. The de~Rham
derivative preserves the above spaces, 
and the quotient topological complex
$\left(\overline{{\cal B}_v}/{\cal B}_v,\bar d\right)$ 
is canonically isomorphic to
$\left(\bar 0_1,d_1\right)$. Moreover we have the 
following known result whose proof is
easy and  does not require that $\F$ be Riemannian and $M$ closed. 

\bl[V.~Sergiescu \cite{Sergiescu86}; see also \cite{Alv2}]\label{l:B/Z} 
The cohomology of the quotient complex ${\cal B}_v/{\cal Z}_{v-1}$ is trivial.
\el

We also have the following canonical identities of topological vector spaces:
\ba
H\left({\cal B}_v/{\cal Z}_{v-1},\bar d\right)&\equiv&
d^{-1}({\cal Z}_{v-1})/d({\cal B}_v)\;,\\ 
H\left(\overline{{\cal B}_v}/{\cal Z}_{v-1},\bar d\right)
&\equiv&d^{-1}({\cal Z}_{v-1})/d\left(\overline{{\cal B}_v}\right)\;,\\
H\left(\overline{{\cal B}_v}/{\cal B}_v,\bar d\right)
&\equiv&d^{-1}({\cal B}_v)/d\left(\overline{{\cal B}_v}\right)\;.
\ea
Then
Lemma~\ref{l:B/Z} means 
$d^{-1}({\cal Z}_{v-1})=d({\cal B}_v)$. So the canonical
map
$$H\left(\overline{{\cal B}_v}/{\cal Z}_{v-1},\bar d\right)\ar 
H\left(\overline{{\cal B}_v}/{\cal B}_v,\bar d\right)$$
can be identified with the equality
$$d^{-1}({\cal Z}_{v-1})/d\left(\overline{{\cal B}_v}\right)
=d^{-1}({\cal B}_v)/d\left(\overline{{\cal B}_v}\right)\;,$$
and is thus an isomorphism of topological vector spaces.

Now (\ref{e:leafwise Hodge}) in Theorem~\ref{t:de Rham} implies that
$\widetilde\Pi_{\cdot,v}$ induces an isomorphism of topological 
vector spaces between 
$\overline{{\cal B}_v}/{\cal Z}_{v-1}$ and 
$\widetilde{\cal H}_1^{\cdot,v}$ such that
$\tilde d_1$ corresponds to $\bar d$. Hence 
$\tilde d_1^2=0$ and  
$H\left(\bar 0_1,d_1\right)\cong 
H\left(\widetilde{\cal H}_1,\tilde d_1\right)
\cong\ker\widetilde\D_1$ as topological
vector spaces, yielding $\ker\widetilde\D_1=0$ 
as desired because $\ker\widetilde\D_1$ is Hausdorff while $H\left(\bar
0_1,d_1\right)$ has the trivial topology.

Observe that the duality stated in
Corollary~\ref{c:E2} follows because, when $M$ is oriented, 
the corresponding Hodge
star operator commutes with $\Pi$ and $\D_1$.

\begin{rem}
Let $E(\lambda)$ and $\widetilde E\left(\tilde\lambda\right)$ 
be the eigenspaces
corresponding to eigenvalues $\lambda$ of $\D_1$ and 
$\tilde\lambda$ of
$\widetilde\D_1$. This eigenspaces have bigradings 
induced by the bigradings of ${\cal
H}_1$ and 
$\widetilde{\cal H}_1$. Then, if $M$ is oriented, 
the Hodge star operator also
induces duality $E(\lambda)^{u,v}\cong E(\lambda)^{q-u,p-v}$, and skew
duality $\widetilde E\left(\tilde\lambda\right)^{u,v}\cong\widetilde
E\left(\tilde\lambda\right)^{q-u-1,p-v+1}$.
\end{rem}

\begin{rem}
If we use appropriate zero order modifications of $D_\perp$, then 
we get the Hodge
theoretic version of the results of \cite{Demetrio}, which 
have important implications
about tenseness.  \end{rem}

\end{document}